\documentclass[12pt]{article}
\usepackage{a4wide}
\usepackage{amsmath,amsfonts,amssymb}
\usepackage{amsthm}
\usepackage{graphicx}

\usepackage{xcolor}

\usepackage{hyperref}
\usepackage{breakurl}
\usepackage[T1]{fontenc}

\usepackage[numbers,sort&compress]{natbib}

\newcommand{\conv}[1]{\mathrm{conv}(#1)}

\newtheorem{theorem}{Theorem}        
            
\newtheorem{lemma}[theorem]{Lemma}

\theoremstyle{remark}

\begin{document}
	\title{Existence of alternating Hamiltonian paths on specific bicolored point sets
	}

	\author{
		Jan Soukup\footnote{Department of Applied Mathematics, Charles University, Faculty of Mathematics and Physics, Malostransk\'e n\'am.~25, 118 00~ Praha 1, Czech Republic; \texttt{soukup@kam.mff.cuni.cz}. Supported by project 23-04949X of the Czech Science Foundation (GA\v{C}R) and by the grant SVV--2025--260822.  
		}
	} 
	\maketitle              
	\begin{abstract}
		Consider a bicolored point set $P$ in general position in the plane consisting of $n$ blue and $n$ red points. We show that if $P$ can be partitioned into bicolored pairs $(a_i,b_i)$ such that, for every $i$, $a_i$ and $b_i$ have different colors and the sets $P\setminus \{a_i\}$ and $P \setminus \{b_i\}$ have the same order type, then there exists a noncrossing alternating Hamiltonian path on $P$. We also investigate the case when $P$ lies on the circle, and every interval along the circle contains at most~$c$ more points of one color than of the other. For $c=2$, we show that there exists a noncrossing alternating Hamiltonian path on $P$. For $c=3$, we show that there exist configurations of points for which the longest noncrossing alternating path on $P$ has size at most $(2-1/25)n + O(1)$.
		
	\end{abstract}

	\section{Introduction} 
	In geometric graph theory, it is a common problem to decide whether a given graph can be drawn in the plane on a given point set so that the edges are represented by non-crossing straight-line segments. For example, deciding whether a given general planar graph admits a~non-crossing straight-line drawing on a given point set is NP-complete~\cite{c06}. 
	
	There are many interesting unanswered questions when considering bicolored point sets instead; see the comprehensive survey by Kano and Urrutia~\cite{ku21}. We focus on drawings of bipartite graphs on bicolored point sets, where edges are drawn as non-crossing straight-line segments between points of different colors. This question remains interesting even for paths. Let $B$ and $R$ denote sets of blue and red points in the plane, respectively, such that $R \cup B$ is in general position, i.e., no three points are collinear. We call a non-intersecting path on $R\cup B$ whose edges are straight-line segments, each connecting two points of $R\cup B$ of different colors, an $\emph{alternating path}$. If such an alternating path connects all points of $R\cup B$, we call it an \emph{alternating Hamiltonian path}.
	
	If $||R| - |B|| \le 1$ and $R$ can be separated from $B$ by a line, then Abellanas et al.~\cite{aghnr99} showed that there always exists an alternating Hamiltonian path on $R \cup B$. This fact, together with the well-known Ham sandwich theorem, implies that if $|R| = |B|$, then there always exists an alternating path on $R \cup B$ connecting at least half of the points. This trivial lower bound on the length of an alternating path that always exists is the best known to our knowledge. This bound was improved by a small linear factor by Mulzer and Valtr \cite{mv20} for point sets in \emph{convex} position, i.e., when the points form the vertices of a convex polygon. On the other hand, if we do not assume that $R$ and $B$ are separated by a line, then there are examples with $|R| = |B| \ge 8$ in which no alternating Hamiltonian path on $R \cup B$ exists, even if $R\cup B$ is in convex position. Moreover, for $R\cup B$ in convex position with $|R| = |B| = n$, Cs\'{o}ka et al. \cite{cbkl22} showed that there are configurations where the longest alternating path on $R\cup B$ has size at most $(4-2\sqrt{2})n  + o(n)$.	
	
    As noted above, an alternating Hamiltonian path does not exist on every point set, but we have already given one example of a sufficient condition on the point set for which an alternating Hamiltonian path exists. Several other sufficient conditions are known~\cite{ckmsv13, aghnr99, s24}. As our first result, we extend the known family of point sets for which there exists an alternating Hamiltonian path on $R \cup B$. Specifically, we prove that it is sufficient if the points of $R\cup B$ can be partitioned into bicolored pairs.

    Two distinct points $x, y$ in a point set $P$ are called \emph{$P$-twins} or \emph{$P$-twin points} if, for every $a,b \in P\setminus\{x,y\}$, the triangle $aby$ has the same orientation as the triangle $abx$.

    \begin{theorem}\label{thm_doubled}
        Let $R = \{r_1, \dots, r_n\}$ be a set of $n$ red points and $B = \{b_1, \dots, b_n\}$ be a set of $n$ blue points such that $R \cup B$ is in general position and, for every $i$, $r_i$ and $b_i$ are $(R \cup B)$-twins. Then $R\cup B$ admits an alternating Hamiltonian path. 
    \end{theorem}

    Most of the bounds on the lengths of alternating paths on $R\cup B$ were not proved directly. Instead, various authors used a connection with so-called separated matchings. 
    
    A matching on $R\cup B$ whose edges are straight-line segments, each connecting two points of $R\cup B$ of different colors, is called an $\emph{alternating matching}$. An alternating matching is called a \emph{separated matching} if there exists a line that intersects the interior of every segment in the matching.
    
    Let $\mu(n)$ and $l(n)$ be the largest numbers such that, for any set $R \cup B$ in convex position with $|R| = |B| = n$, there always exists a separated matching or an alternating path, respectively, covering at least $\mu(n)$ or $l(n)$ vertices, respectively. Again, we may assume that the points lie on the circle. Trivially, the edges of a separated matching can be joined by new edges to form an alternating path. Thus, $\mu(n) \le  l(n)$. The best lower bound on $l(n)$ by Mulzer and Valtr \cite{mv20} was actually proved for separated matchings, so the bound is the same for $\mu(n)$ and  $l(n)$. Kyn\u{c}l, Pach and T\'{o}th \cite{kpt08} showed that a variation of the opposite inequality also holds; they proved that if $k$ is the number of color changes along the circle, then $l(n) -4k -1 \le \mu(n)$. Thus, if the number of such color changes is in $o(n)$, the asymptotics of $l(n)$ and $\mu(n)$ differ only by a function in $o(n)$. All the upper bounds on $l(n)$ mentioned above rely on the same fact; they prove an upper bound on $\mu(n)$ using constructions with $o(n)$ color changes, and so the same upper bound holds for $l(n)$ up to an $o(n)$ difference. 

    To better understand the behavior of $\mu(n)$ independently of $l(n)$, M\'esz\'aros~\cite{m11} studied cases in which the number of color changes along the circle is linear in $n$. In particular, they considered cases with small discrepancies. 
    The \emph{discrepancy} of a set of red and blue points $R \cup B$, $|R| = |B| = n$, on a circle is $d$ if, for every interval on the circle, the difference between the cardinalities of its color classes is at most $d$. M\'esz\'aros~\cite{m11} proved that if the discrepancy of $R\cup B$ is two or three, then $R\cup B$ admits a separated matching of size at least $4/3n$; this is larger than the size that can be achieved in the general case. Soukup~\cite{s21} later improved the case of discrepancy two by showing that $R\cup B$ admits a separated matching of size $8/5n - o(n)$. 

    We continue this study of the distinct asymptotics of $\mu(n)$ and $l(n)$ in the case of small discrepancies. Our main results are the following theorems that bound the parameter $l(n)$ for the cases of discrepancy two and three. 

    \begin{theorem}\label{thm_discrepacy2}
    For every $n$, every point set $S$ consisting of $n$ red and $n$ blue points on the circle with discrepancy two admits a Hamiltonian alternating path. 
    \end{theorem}
    
    \begin{theorem}\label{thm_discrepacy3}
    For every $n$ there exists a set $S$ of $n$ red and $n$ blue points on the circle with discrepancy three such that every alternating path on $S$ has size at most $(2 -\frac{1}{25})n + O(1)$.
    \end{theorem}

    Theorem~\ref{thm_discrepacy3} is the first theorem to prove an upper bound on $l(n)$ using a construction in which the largest monochromatic consecutive subset of points on the circle has constant size. We hope it will help prove a better upper bound on $l(n)$ in the general case.

We prove Theorem~\ref{thm_doubled} in Section~\ref{sec:doubled}, and we prove Theorems~\ref{thm_discrepacy2} and~\ref{thm_discrepacy3} in Section~\ref{sec:discrepancy3}.

\section{Hamiltonian paths on point sets formed by twin points}\label{sec:doubled}
        
    \begin{proof}[Proof of Theorem~\ref{thm_doubled}]
        In this proof, we work with the point set $R\cup B$, so we call every $(R\cup B)$-twin simply a \emph{twin} for brevity. We show how to inductively build a sequence $p_1q_1p_2q_2\dots p_nq_n$ that contains all points of $R\cup B$ such that connecting them in this order with straight-line segments yields an alternating Hamiltonian path. For each $p_i$, let $S_i$ be the set of unselected points so far, that is, $S_i = \{q_i, p_{i+1}, q_{i+1}, \dots, p_n, q_{n}\}$. Furthermore, for each $p_i$, we need an auxiliary point $t_i\in S_i$.
        
        We inductively construct the sequence and auxiliary points so that the following conditions are satisfied for all $i$ with $1\le i \le n$:
        \begin{enumerate}
            \item[\rm 1)] The points $p_i$ and $t_i$ have different colors. 
            \item[\rm 2)] The point $p_i$ can be separated from $S_i$ by a line $l$ such that $l \cap \conv{\{p_i, t_i\}} = \{t_i\}$ and $l \cap \conv{S_i} = \{t_i\}$. 
            \item[\rm 3)] Every point in $S_i$ has its twin in $S_i$, except for $t_i$, whose twin is not in $S_i$.
            \item[\rm 4)] The path $p_1q_1\dots p_i$ is non-crossing and intersects $\conv{S_i \cup\{p_i\}}$ only at $p_i$.
        \end{enumerate}
        Clearly, by the last condition for $i=n$, the path $p_1q_1\dots p_{n-1}q_{n-1}p_nt_n$ is a non-crossing Hamiltonian path on $R \cup B$.

        \begin{figure}
			\centering
			\includegraphics[scale=1]{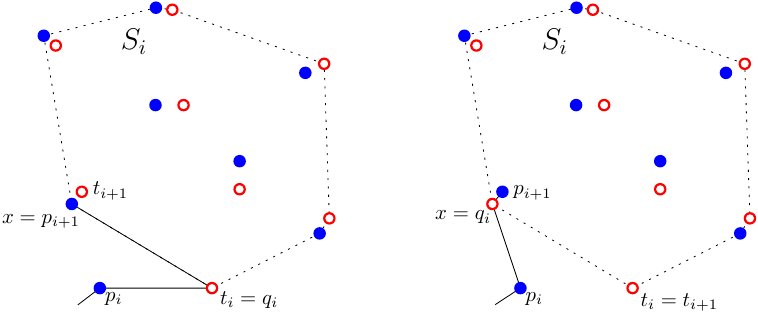}
        \caption[The induction step from Theorem~\ref{thm_doubled}.]{Illustration of how to inductively find an alternating Hamiltonian path in the proof of Theorem~\ref{thm_doubled}. The picture illustrates the induction step. There are two main cases: on the left, point $x$ has the same color as point $p_i$; on the right, they have different colors.}
			\label{fig_doubled}
		\end{figure}

        We now prove that we can construct the sequence satisfying the given conditions. 
        
        We set $p_1$ to be a point on the boundary of $\conv{R\cup B}$ and $t_i$ to be the twin point of~$p_i$. The conditions are clearly satisfied for $i=1$.

		Assume we have already built the sequence up to $p_i$, selected $t_i$, and that the conditions hold for all $j\le i$. If $i=n$, we set $q_n = t_n$ and are done, as discussed earlier. Thus, we can assume $1\le i<n$. Assume $p_i$ is blue; otherwise, we can switch the colors of all points. See Figure~\ref{fig_doubled} for an illustration. By condition~1), $t_i$ is red. By condition~2), there exists a line~$l$ through $t_i$ that separates $p_i$ from $S_i$. By rotating $l$ around $t_i$, we can find a point $x \in S_i$, a vertex of $\conv{S_i}$ distinct from $t_i$, such that the line $t_ix$ separates $S_i$ from $p_i$. Note that the segment $xt_i$ is in fact an edge of $\conv{S_i}$. We split the analysis according to the color of $x$. 

        First, assume that $x$ is blue. We set $q_i = t_i$ and $p_{i+1} = x$. Additionally, we set $t_{i+1}$ to be the twin of $x$ (it is in $S_i$ by condition~3)). See Figure~\ref{fig_doubled} (left) for an illustration. We now verify that all conditions hold for $i+1$. Since $p_{i+1}$ and $t_{i+1}$ are twins, they have distinct colors. Since the segment $t_ix$ is an edge of $\conv{S_i}$, the point $t_{i+1}$ is in $S_i$, and $x$ is a twin of $t_{i+1}$, it follows that the line $t_it_{i+1}$ separates $S_{i+1}$ from the point $p_{i+1}$ and intersects both $\conv{S_{i+1}}$ and the segment $p_{i+1}t_{i+1}$ only at $t_{i+1}$. Thus, condition 2) is satisfied. Since $S_{i+1} = S_i \setminus \{t_i, p_{i+1}\}$ and $p_{i+1}$ and $t_{i+1}$ are twins, condition 3) for $i+1$ follows directly from condition 3) for $i$. Finally, condition 4) holds because the newly added part of the path $p_1q_1\dots p_{i+1}$, namely $p_it_ip_{i+1}$, lies in the previously uncrossed triangle $p_it_ip_{i+1}$, and this triangle intersects $\conv{S_{i+1} \cup \{p_{i+1}\}}$ exactly in $\{p_{i+1}\}$.

        Now assume that $x$ is red. We set $q_i = x$ and $t_{i+1} = t_i$. Additionally, we set $p_{i+1}$ to be the twin of $x$ (it is in $S_i$ by condition~3)). See Figure~\ref{fig_doubled} right for an illustration. We now verify that all conditions hold for $i+1$. The points $p_{i+1}$ and $t_{i+1}$ have distinct colors since $t_{i+1} = t_i$ and $p_i$ and $p_{i+1}$ have the same color. Since the segment $t_ix$ is an edge of $\conv{S_i}$, and the point $p_{i+1}$ is in $S_i$ and is the twin of $x$, it follows that the segment $t_ip_{i+1}$ is an edge of $\conv{S_{i+1} \cup \{p_{i+1}\}}$. By rotating the line $t_{i+1}p_{i+1}$ around $t_{i+1}$ slightly, we can find a line that separates $S_{i+1}$ from the point $p_{i+1}$ and intersects both $\conv{S_{i+1}}$ and the segment $p_{i+1}t_{i+1}$ only at $t_{i+1}$. Thus, condition 2) is satisfied. Since $S_{i+1} = S_i \setminus \{x, p_{i+1}\}$ and $p_{i+1}$ and $x$ are twins, condition 3) for $i+1$ follows directly from condition 3) for~$i$. We now verify condition 4). The newly added part of the path $p_1q_1\dots p_{i+1}$, namely $p_ixp_{i+1}$, lies within the previously uncrossed triangle $p_ixp_{i+1}$. Since the line $t_ix$ separates $S_i$ from the point $p_i$, and $p_{i+1}$ is in $S_i$ and is the twin of $x$, it follows that the line $t_ip_{i+1}$ separates $\conv{S_{i+1} \cup \{p_{i+1}\}}$ from the points $x$ and $p_{i}$. Thus, the line $t_ip_{i+1}$ separates the path~$p_ixp_{i+1}$ from $\conv{S_{i+1} \cup \{p_{i+1}\}}$ and intersects both only at $p_{i+1}$.  
    \end{proof}

\section{Alternating paths on convex sets with small discrepancies}\label{sec:discrepancy3}

The \emph{discrepancy} of a sequence $\rho$ of red and blue points is $d$ if, for every continuous subsequence of $\rho$, the difference between the cardinalities of its color classes is at most $d$. Note that if a bicolored point set $S$ on the circle has discrepancy $d$, and we read $S$ as a sequence in clockwise or counterclockwise order starting at any point, the resulting sequence has discrepancy $d$. 

A maximal monochromatic continuous subsequence in a bicolored sequence is called a \emph{run}. Similarly, a \emph{run} of a bicolored point set on the circle is a maximal monochromatic consecutive subset of points of the point set.

We begin with the proof of Theorem~\ref{thm_discrepacy2}. It is almost a direct consequence of Theorem~\ref{thm_doubled}. We need only one additional observation. The following lemma was already proved by Soukup~\cite{s21} in a slightly different form, but we include the proof for completeness. 

\begin{lemma}\label{lemma_pairsplitoncircle}
    Let $S$ be a set of $n$ red and $n$ blue points on a circle with discrepancy two. Then $S$ can be partitioned into pairs of consecutive points, each containing one red and one blue point. 
\end{lemma}
\begin{proof}
    Let $p_1, p_2, \dots, p_{2n}$ be the points of $S$ on the circle, listed in clockwise order. We will use indices modulo $2n$ for the points on the circle.
    
    Since the discrepancy of $S$ is two, every run in $S$ has size at most 2. 
    If each run has size one, we can go in the clockwise direction and always pair a blue point with its successor, which is red. 
    
    Otherwise, there exist two consecutive points of the same color. By re-indexing, we can assume that these points are $p_1$ and $p_{2}$. We can also assume they are blue. We now split $S$ into pairs $(p_{2}, p_{3}), (p_{4}, p_{5}), \dots, (p_{2n}, p_{1})$ of consecutive points along the circle. Note that each pair starts with a point with an even index. We claim that each pair now contains two points of different colors. For contradiction, assume some pair $(p_j, p_{j+1})$ contains points of the same color. We know that $j$ is even and not equal to $2$ or $2n$. Consider the sequence $p_3p_4\dots p_{j-1}$. It has odd length, and because the discrepancy of $S$ is at most two, $p_3p_4\dots p_{j-1}$ contains either one more blue point than red points or one more red point than blue points. But the first case is not possible; otherwise $p_1p_2p_3p_4\dots p_{j-1}$ would contain three more blue points than red points. Thus, $p_3p_4\dots p_{j-1}$ contains one more red point than blue points. Hence, if $p_j, p_{j+1}$ are red, then $p_3p_4\dots p_jp_{j+1}$ contains three more red points than blue points, a contradiction. Similarly, $p_j, p_{j+1}$ cannot be blue because otherwise $p_1p_2p_3p_4\dots p_jp_{j+1}$ would contain three more blue points than red points. Therefore, no such monochromatic pair $(p_j, p_{j+1})$ can exist.
\end{proof}

\begin{proof}[Proof of Theorem~\ref{thm_discrepacy2}]
    By Lemma~\ref{lemma_pairsplitoncircle}, we can partition the points of $R\cup B$ into pairs of consecutive points, each containing one red and one blue point. The points in each pair are twins in $R \cup B$, so the theorem follows directly from Theorem~\ref{thm_doubled}.
\end{proof}

We continue with the proof of Theorem~\ref{thm_discrepacy3}.
In the proof, we use the following sequence of points. Let $A$ be a sequence of 20 points consisting of eight runs, with the first run blue and runs of sizes $3,2,2,3,3,3,2,2$ from first to last, respectively. 

Additionally, we first prove a lemma about an alternating path on two sequences of points. Let $S_1$ be a sequence of red and blue points on a line, and let $S_2$ be a sequence of red and blue points on another line parallel to the first one. Let $P$ be an alternating path on the points of $S_1$ and $S_2$ such that the points of $S_1$ on $P$ appear in the same order as in $S_1$, and the same holds for the points of $S_2$. Then we call $P$ \emph{dense} if the points of $P$ additionally form a continuous subsequence in both $S_1$ and $S_2$. Furthermore, a sequence of points on a line is \emph{double-infinite} if it is infinite in both directions.

\begin{lemma}\label{lemma_twointervalsdis3}
    Let $S_1$ be a double-infinite sequence of points on a line repeating the sequence $A$, and let $S_2$ be a double-infinite sequence of points repeating the sequence $A^T$ on another line that is parallel to the first line. Then the longest alternating dense path between $S_1$ and $S_2$ has size $48$.
\end{lemma}

\begin{proof}
    Let $\dots, A_{-1}, A_0, A_1, \dots $ be the runs of $S_1$ such that $S_1 = \dots A_{-1}A_0A_1 \dots $ and $A_1\dots A_8 = A$. Similarly, let $\dots, B_{-1}, B_0, B_1, \dots $ be the runs of $S_2$ such that $S_2 = \dots B_{-1}B_0B_1 \dots $ and $B_1\dots B_8 = A^T$. Let $\rho$ be a dense alternating path between $S_1$ and $S_2$. 

    First, we observe that if at least two points of a run $X$ lie in $\rho$, then some point in $X$ must be connected by an edge to a run from the other sequence. Let $x_1$ and $x_2$ be two points of $X$ that appear on $\rho$ in this order, and let $y$ be the point immediately after $x_1$ in $\rho$ (it exists because $\rho$ is dense and contains $x_2$). Since $\rho$ is dense, either $y=x_2$ or $y$ is not from the same sequence as $x_1$. But $x_2$ has the same color as $x_1$, which differs from the color of $y$. Thus, $y$ is not from the same sequence as $x_1$. 
    
    Now, we observe that every run of $S_1$ can be connected in $\rho$ by edges to at most one run in $S_2$. Suppose that is not the case, so a run $X$ in $S_1$ is connected to two runs $Y_1, Y_2$ of $S_2$ by edges $e_1$ and $e_2$. The runs $Y_1, Y_2$ cannot have the same color as $X$. Thus, there is a run $Y_3$ of the same color as $X$ between $Y_1$ and $Y_2$. Since $\rho$ is dense, all points of $Y_3$ are in $\rho$, so $Y_3$ must be connected to $S_1$ by an edge $e_3$. But that is not possible because $e_3$ cannot cross $e_1$ and $e_2$, and $e_3$ cannot be adjacent to $X$ because $X$ and $Y_3$ have the same color. Similarly, we can prove that every run of $S_2$ can be connected in $\rho$ by edges to at most one run in $S_1$.

    \begin{figure}
			\centering
			\includegraphics[scale=1]{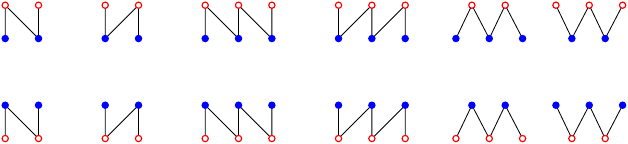}
        \caption[Possible matching pairs of runs from the proof of Lemma~\ref{lemma_twointervalsdis3}.]{All possible matching pairs of runs from the proof of Lemma~\ref{lemma_twointervalsdis3}, along with all possible paths connecting them.}
			\label{fig_matchingpairs}
	\end{figure}
        
    Therefore, if at least two points of a run $X$ of $S_1$ or $S_2$ lie in $\rho$, then $X$ is connected to exactly one run $Y$ of the other $S_i$. 
    If all points of the runs $X$ and $Y$ are in $\rho$, we say that $X$ and $Y$ are \emph{matched}. Clearly, a run can be matched only to runs of a different color. Furthermore, since $\rho$ is dense, if $A_i$ is matched to $B_j$, then for every $k>1$, the run $A_{i+k}$ can be matched only to $B_{j+k}$. Additionally, since $\rho$ is dense, the restriction of $\rho$ to the matching pair is a subpath of $\rho$. Figure~\ref{fig_matchingpairs} shows all possible ways in which pairs of runs can be matched together. 

    \begin{figure}
			\centering
			\includegraphics[scale=1]{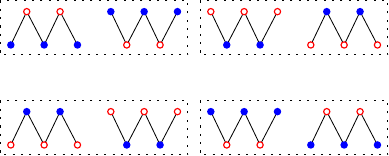}
        \caption[All possible forbidden $2\text{-}3$ configurations from Lemma~\ref{lemma_twointervalsdis3}.]{All possible forbidden $2\text{-}3$ configurations from Lemma~\ref{lemma_twointervalsdis3}. The matching paths cannot be connected together because their endpoints have the same colors.}
			\label{fig_forbiddenpairs}
	\end{figure}
    
    It follows that a pair of runs $A_iA_{i+1}$ of sizes $2$ and $3$, respectively, cannot be matched to a pair of runs $B_jB_{j+1}$ of sizes $3$ and $2$, respectively, because there is only one way to match $A_i$ to $B_{j}$ and $A_{i+1}$ to $B_{j+1}$, and there is no way to connect the matching paths together. Similarly, a pair of runs $A_iA_{i+1}$ of sizes $3$ and $2$, respectively, cannot be matched to a pair of runs $B_jB_{j+1}$ of sizes $2$ and $3$, respectively. See Figure~\ref{fig_forbiddenpairs}. We will call these configurations \emph{forbidden $2\text{-}3$ configurations}.

    \begin{figure}
			\centering
			\includegraphics[scale=1]{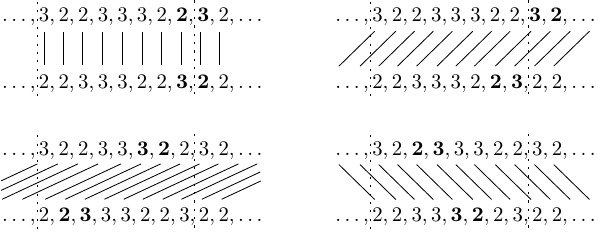}
        \caption[All possible matchings between runs of two sequences.]{All possible matchings between runs of two sequences. The top sequence is the double-infinite sequence repeating $A$, and the bottom one is the double-infinite sequence repeating $A^T$. The sequences are represented by the lengths of their runs. Their periods are separated by dotted lines.}
			\label{fig_periodmatching}
	\end{figure}
    Now suppose that a continuous subsequence of $S_1$ formed by runs of $S_1$ is matched to a continuous subsequence of $S_2$ formed by runs of $S_2$. Since both $S_1$ and $S_2$ have a run period of eight and the colors of runs in both sequences alternate, there are only four ways the matching can look with respect to the period; see Figure~\ref{fig_periodmatching}. In all cases, a forbidden $2\text{-}3$ configuration occurs with a run period of eight. Therefore, $\rho$ contains at most eight matching pairs of runs.

    We can now complete the proof of the lemma. First, assume that $\rho$ contains at most seven matching pairs of runs. These pairs account for at most 36 points. Consider the subpath $\rho_1$ of $\rho$ that precedes these matching pairs. Since $\rho_1$ contains no matching pairs and is dense, it cannot contain an entire run from one of the sequences, say $S_1$. Thus, $\rho_1$ contains at most two points from $S_1$. Additionally, since every run of $S_2$ with at least two points in $\rho_1$ is connected to a run of $S_1$, it follows that $\rho_1$ contains at most three points from a single run of $S_2$ and at most one additional point of $S_2$ from another run of $S_2$. In total, $\rho_1$ has at most six points. Symmetrically, the part of $\rho$ after the matching pairs contains at most six points. Thus, $\rho$ has at most 48 points.    

    \begin{figure}
			\centering
			\includegraphics[scale=1]{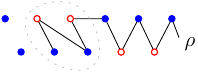}
        \caption[The start of the path $\rho$ from Lemma~\ref{lemma_twointervalsdis3}.]{The start of path $\rho$ from Lemma~\ref{lemma_twointervalsdis3}. If $\rho$ has eight matching pairs of runs, the first pair is a part of a forbidden $2\text{-}3$ configuration. In the picture, we have a forbidden $2\text{-}3$ configuration, and it is easy to see that the highlighted part of $\rho$ before the first matching pair of runs can contain at most four points.}
			\label{fig_pathstart}
	\end{figure}
    
    Now, assume that $\rho$ contains exactly eight matching pairs of runs. Since $\rho$ cannot contain a forbidden $2\text{-}3$ configuration, but a forbidden $2\text{-}3$ matching occurs with a run period of eight, it means that the first matching pair of $\rho$ must be the last matching pair of a forbidden $2\text{-}3$ configuration, and the last matching pair of $\rho$ must be the first matching pair of a forbidden $2\text{-}3$ configuration. All cases are analogous up to recoloring and symmetry, and since $\rho$ is dense, it is easy to see that $\rho$ contains at most four points before the matching pairs and at most four points after them; see Figure~\ref{fig_pathstart} for an illustration.       
\end{proof}

We are ready to prove Theorem~\ref{thm_discrepacy3}.

\begin{proof}[Proof of Theorem~\ref{thm_discrepacy3}]
    Let $n' = \lfloor \frac{n}{20} \rfloor$ and $t = n \bmod 20$. Let $S'$ be the sequence $(br)^tA^{2n'}$, where $b$ denotes a blue point and $r$ denotes a red point. We call the first $2t$ points \emph{exception points}. We place $S'$ in clockwise order on a circle and denote the resulting point set $S$. Note that $S$ has exactly $n$ blue points and $n$ red points. We prove that $S$ satisfies the theorem.   
    
    Let $P = p_1\dots p_m$ be the longest alternating path on $S$. Let $l$ be the line that crosses the first and last edge of $P$. Let $S_1$ be the sequence of points in $S$ on one arc determined by $l$ in clockwise order, and let $S_2$ be the sequence of points in $S$ on the other arc determined by $l$ in counterclockwise order. Then $P$ is an alternating path on the points of $S_1$ and $S_2$ such that the points of $S_1$ on $P$ appear in the same order as in $S_1$, and the same holds for points of $S_2$. Furthermore, we can straighten the two arcs into parallel lines, and the path remains non-crossing. Note that $S_1$ and $S_2$, excluding the exception points, are intervals of the double-infinite sequences repeating $A$ and $A^T$, respectively. We split $P$ into $\lfloor m / 49\rfloor$ subpaths of length at least 49. Lemma~\ref{lemma_twointervalsdis3} implies that no subpath of $P$ of size $49$ without exception points is dense between $S_1$ and $S_2$. Since there are $2t$ exception points, at least $\lfloor m / 49\rfloor -2t$ subpaths are not dense between $S_1$ and $S_2$. For each subpath of $P$ that is not dense, there exists a unique point of $S$ that is not on $P$. Therefore, $2n \ge m + \lfloor m / 49\rfloor -2t$. Thus, $m\le (2-1/25)n + O(1)$. 
\end{proof}

	\bibliographystyle{splncs04nat}
	\bibliography{mybibliography}

\end{document}